\documentclass[a4paper,12pt]{amsart}
\usepackage{amsmath}
\usepackage{amsfonts}
\usepackage{amssymb}
\usepackage{hyperref}
\usepackage{fourier}
\usepackage{a4wide}
\usepackage[ddmmyyyy]{datetime}
\usepackage{tikz}

\newtheorem{lemma}{Lemma}
\newtheorem{theorem}[lemma]{Theorem}
\newtheorem{corollary}[lemma]{Corollary}

\newcounter{rotcount}
\newtheorem{claim}{}[rotcount]

\newenvironment{proofof}{\noindent}{\hfill$\Box$\medskip}
\newenvironment{claimproof}{\noindent {\it Subproof: }}{\hfill$\lozenge$\smallskip}

\newcommand{\ncont}{\nsubseteq}\newcommand{\cont}{\subseteq}
\renewcommand{\u}{\cup}
\newcommand{\s}{^*}\newcommand{\del}{\backslash}
\newcommand{\cl}{{\rm cl}}
\newcommand{\si}{{\rm si}}\newcommand{\co}{{\rm co}}

\newcommand{\defin}{\textbf}

\title{Triangle-roundedness in matroids}
\author[Jo\~ao Paulo Costalonga]{Jo\~ao Paulo Costalonga$^1$}
\thanks{$^1$Departamento de Matem\'atica, Universidade Federal do Esp\'irito Santo. Av. Fernando Ferrari, 514; Campus de Goiabeiras,
  29075-910, Vit\'oria, ES, Brazil. e-mail: {\upshape joaocostalonga@gmail.com} (corresponding author).}
\author[Xiangqian Zhou]{Xiangqian Zhou$^2$}
\thanks{$^2$Department of Mathematics and Statistics, Wright State University, Dayton, OH, 45435, USA and School of Mathematical Sciences, Huaqiao University, Fujian, China. e-mail: {\upshape xiangqian.zhou@wright.edu}}
 
 \newcommand{\M}{\mathcal{M}}
\newcommand{\F}{\mathcal{F}}
\newcommand{\W}{\mathcal{W}}

\begin{document}
\begin{abstract}
A matroid $N$ is said to be triangle-rounded in a class of matroids $\mathcal{M}$ if each $3$-connected matroid $M\in \mathcal{M}$ with a triangle $T$ and an $N$-minor has an $N$-minor with $T$ as triangle. Reid gave a result useful to identify such matroids as stated next: suppose that $M$ is a binary $3$-connected matroid with a $3$-connected minor $N$, $T$ is a triangle of $M$ and $e\in T\cap E(N)$; then $M$ has a $3$-connected minor $M'$ with an $N$-minor such that $T$ is a triangle of $M'$ and $|E(M')|\le |E(N)|+2$. We strengthen this result by dropping the condition that such element $e$ exists 

and proving that there is a $3$-connected minor $M'$ of $M$ with an $N$-minor $N'$ such that $T$ is a triangle of $M'$ and $E(M')-E(N')\subseteq T$. This result is extended to the non-binary case and, as an application, we prove that $M(K_5)$ is triangle-rounded in the class of the regular matroids.
\end{abstract}

\maketitle

\section{Introduction}

Let $\M$ be a class of matroids closed for minors and isomorphisms and let $\F$ be a family of matroids. An \defin{$\F$-minor} of a matroid $M$ is a minor of $M$ isomorphic to a member of $\F$. A matroid $M$ \defin{uses} a set $T$ if $T\cont E(M)$. We say that $\F$ is \defin{$(k,t)$-rounded} in $\M$ if each element of $\F$ is $k$-connected and, for each $k$-connected matroid $M\in\M$ with an $\F$-minor and each $t$-subset $T\cont E(M)$, $M$ has an $\F$-minor using $T$. We define $\F$ to be \defin{$t$-rounded} in $\M$ if it is $(t+1,t)$-rounded in $\M$. A matroid $N$ is said $(k,t)$-rounded (resp. $t$-rounded) in $\M$ if so is $\{N\}$. When we simply say that a matroid or family of matroids is $(k,t)$-rounded or $t$-rounded with no mention to a specific class of matroids, we are referring to the class of all matroids.

Bixby~\cite{Bixby} proved that $U_{2,4}$ is $1$-rounded. Seymour~\cite{Seymour1977} established a method to find a minimal $1$-rounded family containing a given family of matroids; in that work it is established that $\{U_{2,4}$,$M(K_4)\}$, $\{U_{2,4}$,$F_7$,$F_7^*\}$, $\{U_{2,4}$,$F_7$,$F_7^*$,$M^*(K_{3,3})$,$M^*(K_5)$,$M^*(K_{3,3}')\}$ and $\{U_{2,5},U_{3,5},F_7,F_7^*\}$ are $1$-rounded.

Seymour~\cite{Seymour1981} proved that $U_{2,4}$ is also $2$-rounded and, later, in \cite{Seymour1985}, established a method to find a minimal $2$-rounded family containing a given family of matroids.

Khan~\cite{Khan1985} and Coullard~\cite{Coullard1986} proved independently that $U_{2,4}$ is not $3$-rounded. To the best of our knowledge, there is no known criterion to check $(k,t)$-roundedness for $k\ge 4$. For $t\ge 3$, Oxley~\cite{Oley1987} proved that $\{U_{2,4},\mathcal{W}^3\}$ is $(3,3)$-rounded. Moss~\cite{Moss} proved that $\{\mathcal{W}^2, \W^3, \W^4, M(\W_3), M(\W_4),Q_6\}$ is $(3,4)$-rounded and $\{M(\W_3)$,$M(\W_4)$,$M(\W_5)$, $M(K_5\del e)$,$M\s(K_5\del e)$,$M(K_{1,2,3})$,$M\s(K_{1,2,3})$,$S_8\}$ is $(3,5)$-rounded in the class of the binary matroids.

There are results on classification of small $t$-rounded families of matroids for $t=1,2$. Oxley~\cite{Oxley1984} proved that for $|E(N)|\ge4$, $N$ is $1$-rounded if and only if $N\cong U_{2,4}, P(U_{1,3},U_{1,1})$ or $Q_6$ and $2$-rounded if and only if $N\cong U_{2,4}$. Reid and Oxley~\cite{OxleyReid1990} proved that, up to isomorphisms, the unique $2$-rounded matroids with more than three members in the class of $GF(q)$-representable matroids are $M(\mathcal{W}_3)$ and $\M(\mathcal{W}_4)$ for $q=2$, $U_{2,4}$ and $\mathcal{W}^3$ for $q=3$ and $U_{2,4}$ for $q\ge4$.

In this work, we focus on a different type of roundedness. A family of matroids $\F$ is said to be \defin{triangle-rounded} in $\M$ if all members of $\F$ are $3$-connected and, for each matroid $M\in \M$ with an $\F$-minor and each triangle $T$ of $M$, there is  an $\F$-minor of $M$ with $T$ as triangle. We say that a matroid $N$ is \defin{triangle-rounded} in $\M$ if so is $\{N\}$. Some examples of triangle-rounded matroids and families are $U_{2,4}$ in the class of all matroids, $F_7$ in the class of binary matroids and $M\s(K_{3,3})$ in the class of regular matroids (Asano, Nishizeki and Seymour~\cite{Asano}), $M(K_5\del e)$ in the class of regular matroids and $\{S_8,J_{10}\}$ in the class of binary matroids (Reid~\cite{Reid}). The proofs for the triangle-roundedness of the later two rely on the following criterion:

\begin{theorem}(Reid \cite[Theorem 1.1]{Reid})
Let $\{e, f, g\}$ be a triangle of a $3$-connected binary matroid $M$ and $N$ be a $3$-connected minor of $M$ with $e\in E(N)$. Then, there exists a $3$-connected minor $M'$ of $M$ using $\{e, f, g\}$ such that $M'$ has a minor which is isomorphic to $N$ and $E(M')$ has at most $|E(N)| + 2$ elements. 
\end{theorem}

Here we establish a stronger result for binary matroids:

\begin{theorem}\label{binary-criterion}
If $M$ is a $3$-connected binary matroid with a $3$-connected minor $N$ and $T$ is a triangle of $M$, then $M$ has a $3$-connected minor $M'$ using $T$ with an $N$-minor $N'$ such that $E(M')-E(N')\cont T$.
\end{theorem}

Theorem \ref{binary-criterion} will be proved using a stronger result:

\begin{theorem}\label{criterion}
Let $M$ be a matroid with a $3$-connected minor $N$ satisfying $|E(N)|\ge 4$. Suppose that $T$ is a triangle of $M$ and $M$ is minor-minimal with the property that $M$ is $3$-connected and has an $N$-minor using $T$. Then $r(M)-r(N)\le 2$ and for some $N$-minor $N'$ of $M$, $|E(M)-(E(N')\cup T)|\le 1$. Moreover, if $E(M)-E(N')\ncont T$, then one of the following assertions holds:
\begin{enumerate}
 \item [(a)]$r(M)-r(N)=1$, $r\s(M)-r\s(N)\in\{2,3\}$ and $M$ has an element $x$ such that $E(M)-E(N')\cont T\cup x$, $T\cup x$ is a $4$-point line of $M$, $M\del x$ has no $N$-minor and $x$ is the unique element of $M$ such that $\si(M/x)$ is $3$-connected with an $N$-minor; or
 \item [(b)]$r\s(M)-r\s(N)=2$, $r(M)-r(N)\in\{1,2\}$ and $M$ has an element $y$ such that $E(M)-E(N')\cont T\cup y$, $T\u y$ is a $4$-cocircuit of $M$, $M/y$ has no $N$-minor and $M\del A$ has no $N$-minor for each subset $A$ of $T$ with $|A|\ge 2$.
\end{enumerate}
\end{theorem}

All the possible cases described in this theorem indeed occur, we give examples in Section \ref{sec-sharp}. We say that a graph $G$ is \defin{triangle-rounded} if so is $M(G)$ in the class of graphic matroids. Using Theorem \ref{binary-criterion}, we establish that $K_5$ is triangle-rounded, in other words:

\begin{theorem}\label{k_5 rounded}
If $G$ is a $3$-connected graph with a triangle $T$ and a $K_5$-minor, then $G$ has a $K_5$-minor with $E(T)$ as edge-set of a triangle.
\end{theorem}

\noindent{\bf Remark: } $K_3$ and $K_4$ are triangle-rounded, but no larger complete graph than $K_5$ is triangle-rounded. Indeed, consider, for disjoint sets $X$, $Y$ and $\{z\}$ satisfying $|X|,|Y|\ge 2$, a complete graph $K$ with $n$ vertices with $X\cup Y\cup\{z\}$ as vertex set. Consider also  a graph $G$ extending $K$ by two vertices $x$ and $y$ with $E(G)-E(K)=\{xy,xz,yz\}\cup\{xx':x'\in X\}\cup\{yy':y'\in Y\}$.
Note that $G\del xz/xy\cong K_n$. But no $K_n$-minor of $G$ uses $\{xy,xz,yz\}$ because contracting any other edge than $xy$ in $G$ results in a graph with more than one parallel pair of edges.

\smallskip

The next result allows us to derive triangle-roundedness in the class of regular matroids from triangle-roundedness in the classes of graphic and cographic matroids.

\begin{theorem}\label{regular}
If a family $\F$ of internally $4$-connected matroids with no triads is triangle-rounded both in the class of graphic and cographic matroids, then $\F$ is triangle-rounded in the class of regular matroids not isomorphic to $R_{10}$.
\end{theorem}

As $R_{10}$ has no $M(K_5)$-minor and $M(K_5)$ is internally $4$-connected with no triads and trivially triangle-rounded in the class of cographic matroids, it follows from Theorems \ref{k_5 rounded} and \ref{regular} that:

\begin{corollary}
$M(K_5)$ is triangle-rounded in the class of regular matroids.
\end{corollary}

All proofs are in the next section.

\section{Proofs}

In this section we prove the theorems. For matroids $M$ and $N$ we write $N<M$ to say that $M$ has a proper minor isomorphic to $N$. The notation $N\le M$ means that $N<M$ or $N=M$. Next we state some results used in the proofs.

\begin{lemma}\label{w36}(Whittle, \cite[Lemma 3.6]{Whittle})
Let $M$ be a $3$-connected matroid with elements $x$ and $p$ such that $si(M/x)$ and $si(M/x,p)$ are $3$-connected, but $si(M/p)$ is not $3$-connected. Then, $r(M)\geq 4$ and there is a rank-$3$ cocircuit $C\s$ of $M$ containing $x$ such that $p\in cl_M(C\s)-C\s$.
\end{lemma}

\begin{lemma}\label{w37}(Whittle \cite[Lemma 3.7]{Whittle})
Let $C\s$ be a rank-$3$ cocircuit of a $3$-connected matroid $M$ such that $p\in cl_M(C\s)-C\s$.
\begin{enumerate}
\item[(a)] If $z_1,z_2\in C\s$, then  $si(M/p,z_1)\cong si(M/p,z_2)$.
\item[(b)] If $N$ is a matroid and for some $x\in C\s$, $\si(N/x,p)$ is $3$-connected with an $N$-minor, then $\si(N/z,p)$ is $3$-connected with an $N$-minor for each $z\in C\s$.
\end{enumerate}
\end{lemma}

\begin{lemma}\label{w38}(Whittle \cite[Lemma 3.8]{Whittle})
Let $C\s$ be a rank-$3$ cocircuit of a $3$-connected matroid $M$. If $x\in C\s$ has the property that $cl_M(C\s)-x$ contains a triangle of $M/x$, then $si(M/x)$ is 3-connected.
\end{lemma}

From Lemma \ref{w38}, we may conclude:

\begin{corollary}\label{w38-cor}
Let $M$ be a $3$-connectced matroid with a triangle $T$ and a triad $T\s$ such that $T\s-T=\{x\}$ and $T-T\s=\{y\}$. Then $\si(M/x)$ and $\co(M\del y)$ are $3$-connectced.
\end{corollary}

\begin{lemma}\label{wu}(Wu, \cite[Lemma 3.15]{Wu})
If $I\s$ is a coindependent set in a matroid $M$ and $M\del I\s$ is vertically $3$-connected, then so is $M$.
\end{lemma}

Using Seymour's Splitter Theorem (as stated in \cite[Corollary 12.2.1]{Oxley}) and proceeding by induction on $i$ using Lemma \ref{wu}, we may conclude:

\begin{corollary}\label{splitter-plus}
Let $N<M$ be $3$-connected matroids such that $M$ has no larger wheel or whirl-minor than $N$ in case $N$ is a wheel or whirl respectively. Then, there is a chain of $3$-connected matroids $N\cong M_n<\cdots<M_1<M_0=M$ such that for each $i=1,\dots,n$ there is $x_i\in E(M_i)$ satisfying $M_i=M_{i-1}/x_i$ or $M_i=M_{i-1}\del x_i$. Moreover, for $I:=\{x_i:M_{i-1}=M_i/x_i\}$ and $I\s:=\{x_i:M_{i-1}=M_i/x_i\}$,
\begin{enumerate}
 \item [(a)] $I$ is an independent set and $I\s$ is a coindependent set of $M$.
 \item [(b)] for each $1\le i\le n$, $M/(I\cap\{x_1,\dots,x_i\})$ and $\big(M\del (I\s\cap\{x_1,\dots,x_i\})\big)\s$ are vertically $3$-connected.
\end{enumerate}
\end{corollary}

\begin{theorem}\label{wcomut}(Whittle, \cite[Corollary 3.3]{Whittle}) Let $N$ be a $3$-connected minor of the $3$-connected matroid $M$. If $r(M)\ge r(N)+3$, then for each element $x$ such that $\si(M/x)$ is $3$-connected with an $N$-minor, there exists $y \in E(M)$ such that $si(M/y)$ and $si(M/x, y)$ are $3$-connected with $N$-minors.
\end{theorem}

\begin{theorem}\label{teorema-costalonga}(Whittle \cite[Lemma 3.4 and Theorem 3.1]{Whittle} and Costalonga \cite[Theorem 1.3]{Costalonga2})
Let $k\in\{1,2,3\}$ and let $M$ be a $3$-connected matroid with a $3$-connected minor $N$ such that $r(M)-r(N)\ge k$. Then $M$ has a independent $k$-set $J$ such that $\si(M /x)$ is $3$-connected with an $N$-minor for all $x\in J$.
\end{theorem}

\begin{lemma}(Costalonga \cite[Corollary 4]{Costalonga3})\label{gap4}
Suppose that $N<M$ are $3$-connected matroids with $r\s(M)-r\s(N)\ge 4$ and $N$ is cosimple. Then:
\begin{enumerate}
 \item [(a)] $M$ has a coindependent set $S$ of size $4$ such that $\co(M\del e)$ is $3$-connected with an $N$-minor for all $e\in S$; or
 \item [(b)] $M$ has distinct elements $a_1,a_2,b_1,b_2,b_3$ such that, $T_s:=\{a_s,b_t,b_3\}$ is a triangle for $\{s,t\}=\{1,2\}$, $T\s:=\{b_1,b_2,b_3\}$ is a triad of $M$ and $\co(M\del T\s)$ is $3$-connected with an $N$-minor.
\end{enumerate}
\end{lemma}

Next, we prove Theorem \ref{criterion} and, after, using Theorem \ref{criterion}, we prove Theorem \ref{binary-criterion}. For a non-negative integer $k$, a \defin{$k$-segment} of a matroid is a $k$-subset of a line of this matroid.\\

\begin{proofof}{\it Proof of Theorem \ref{criterion}: }
Suppose that the result does not hold. This is, for each $N$-minor $N'$ of $M$, $E(M)-E(N')\ncont T$ and items (a) and (b) of the theorem do not hold. It is already known that $U_{2,4}$ is triangle-rounded~\cite{Asano}, so, we may assume that $|E(N)|\ge 5$. The proof will be based on a series of assertions. First, note that it follows from the minimality of $M$ that:

\begin{claim}\label{criterion-1}
If, for $x\in E(M)$, $\si(M/x)$ is $3$-connected with an $N$-minor, then $x\in \cl_M(T)$.
\end{claim}

\begin{claim}\label{criterion-2}
If $T\s$ is a triad and $T$ is a triangle of $M$ such that $T\s-T=\{x\}$, then $M\del x$ has no $N$-minor.
\end{claim}
\begin{claimproof}
Suppose the contrary. Let $T\s\cap T=\{a,b\}$. As $N$ is simple and cosimple, $M\del x/a\del b$ has an $N$-minor. But $M\del x/a\del b\cong M\del a,b/x$ and, therefore, $M/x$ has an $N$-minor. By Corollary \ref{w38-cor}, $\si(M/x)$ is $3$-connected. By \ref{criterion-1}, $x\in\cl_M(T)$. As $x\notin T$, then $T\s$ meets a $4$-segment of $M$. This implies that $M\cong U_{2,4}$, a contradiction.
\end{claimproof}

\begin{claim}\label{criterion-3}
If, for $x\in E(M)$, $\co(M\del x)$ is $3$-connected with an $N$-minor, then $x\in T$.
\end{claim}
\begin{claimproof}
Suppose the contrary. Then $T\ncont E(\co(M\del x))$ and therefore, there is a triad $T\s$ meeting $x$ and $T$. This contradicts \ref{criterion-2}.
\end{claimproof}

\begin{claim}\label{criterion-5}
If $\si(M/x)$ and $\si(M/x,y)$ are $3$-connected with $N$-minors, then so is $\si(M/y)$ and $x,y\in T$.
\end{claim}
\begin{claimproof} First we prove that $\si(M/y)$ is $3$-connected. Suppose the contrary. By Lemma \ref{w36}, there is a rank-$3$ cocircuit $C\s$ such that $x\in C\s$ and $y\in\cl(C\s)-C\s$. By \ref{criterion-1}, $x\in \cl(T)$. By orthogonality~\cite[Proposition 2.1.11]{Oxley}, $T\cont\cl(C\s)$. As $r(C\s)=3$, there is $z\in C\s-\cl(T)$ and $T$ is a triangle of $M/z$ contained in $\cl_{M/z}(C\s)$. By Lemma \ref{w38}, $\si(M/z)$ is $3$-connected and, by Lemma \ref{w37}, $M/z$ has an $N$-minor. This contradicts \ref{criterion-1}. So, $\si(M/y)$ is $3$-connected.

By \ref{criterion-1}, $x,y\in \cl(T)$. If, for some $\{a,b\}=\{x,y\}$, $a\notin T$, then, as $M/b$ has an $N$-minor and $a$ is in a parallel pair of $M/b$, it follows that $M\del a$ has an $N$-minor. Moreover, in this case, $T\cup a$ is a $4$-segment of $M$ and $M\del a$ is $3$-connected with an $N$-minor, contradicting the minimality of $M$. Thus, $x,y\in T$.
\end{claimproof}

\begin{claim}\label{criterion-rankgap}
$r(M)-r(N)\le 2$ and $r\s(M)-r\s(N)\le3$.
\end{claim}

\begin{claimproof}
If $r(M)-r(N)\ge 3$, then, by Theorem \ref{teorema-costalonga}, there is an independent set $J$ of size $3$ such that $\si(M/x)$ is $3$-connected with an $N$-minor for all $x\in J$. So, there is $x\in J-\cl_M(T)$, a contradiction to \ref{criterion-1}. Thus, $r(M)-r(N)\le 2$. 

If $r\s(M)-r\s(N)\ge 4$, then Lemma \ref{gap4} applies. If item (a) of that Lemma holds, then we have an element $x\in E(M)-T$ such that $\co(M\del x)$ is $3$-connected with an $N$-minor, contradicting \ref{criterion-3}. So, consider the elements given by item (b) of Lemma \ref{gap4}. Let $s\in\{1,2\}$. Note that $M\del T\s$ is isomorphic to a minor of $M\del a_s$. By Corollary \ref{w38-cor}, $\co(M\del a_s)$ is $3$-connected with an $N$-minor and, therefore, $a_s\in T$ by \ref{criterion-3}. So, $a_1,a_2\in T$. By orthogonality between $T$ and $T\s$, $b_3\notin T$. As $M\del b_3$ has an $N$-minor, then, by \ref{criterion-3}, $\co(M\del b_3)$ is not $3$-connected. By Bixby's Lemma~\cite[Lemma 8.7.3]{Oxley}, $\si(M/b_3)$ is $3$-connected. As $M\del T\s$ has an $N$-minor, then, so has $M/b_3$ and, therefore, $\si(M/b_3)$. By \ref{criterion-1}, $b_3\in \cl_M(T)$. But $b_3$ is in a triad and, therefore, in no $4$-segment of $M$. So, $b_3\in T$, a contradiction.
\end{claimproof}

\begin{claim}\label{criterion-hypothesis}
The theorem holds if there is a wheel or whirl $W$ such that $N<W\le M$.
\end{claim}
\begin{claimproof}
Let $W\cong M(W_n)$ or $W^n$. Since $N$ is a minor of $W$, then $N$ is isomorphic to a wheel or whirl. As $|E(N)|\ge 5$, then $|E(N)|\ge 6$ and $n\ge 4$. As $r(M)-r(N)\le 2$ and $r(W)-r(N)\ge 1$, hence $r(M)-r(W)\le 1$.

Suppose for a contradiction that $r(M)=r(W)$. So, for some coindependet set $J\s$ of $M$, $W=M\del J\s$. If $y$ is a non-spoke of $W$, then $W/y$ is vertically $3$-connected by Corollary \ref{w38-cor}. Moreover, $\si(W/y)$ has an $N$-minor and so does $\si(M/y)$. So, by Lemma \ref{wu}, $\si(M/y)$ is $3$-connected with an $N$-minor for each non-spoke $y$ of $W$. By \ref{criterion-1}, all non-spokes of $W$ are in $\cl_M(T)$, contradicting the fact that they are the elements of a set with rank at least $3$ in $W$.

Thus, $r(M)=r(W)+1$. Now, by Theorem \ref{teorema-costalonga}, there is an element $x$ such that $\si(M/x)$ is $3$-connected with an $W$-minor $W'$. As $r(M/x)=r(W')$, there is coindependent set $I\s$ of $M/x$ such that $W'=M/x\del I\s$. If $y$ is a non-spoke of $W'$, $W'/y=M/x,y\del I\s$ is vertically $3$-connected and, by Lemma \ref{wu}, so is $M/x,y$. By \ref{criterion-5}, all non-spokes of $W'$ are in $T$, a contradiction again. \end{claimproof}

As the theorem fails for $M$, $N$ and $T$, by \ref{criterion-hypothesis}, there is no wheel or whirl $W$ such that $N<W\le M$ and the hypotheses of Seymour's Splitter Theorem hold for $M$ and $N$. 

\begin{claim}\label{criterion-4}
If $x\in E(M)$ and $\si(M/x)$ is $3$-connected with an $N$-minor, then $x\in T$.
\end{claim}
\begin{claimproof}
Suppose the contrary. By \ref{criterion-1}, $x\in \cl_M(T)$, which is a line with more than $3$ points. As $M\del z$ is $3$-connected for all $z\in \cl_M(T)-T$, then $M\del z$ has no $N$-minor if $z\in \cl_M(T)-T$. As $M/x$ has an $N$-minor, $\cl_M(T)=T\cup x$ and $M\del x$ has no $N$-minor. This implies that $r\s(M)-r\s(N)\ge 2$ as $T$ is in a parallel class of $M/x$.

Let us check that for each $z\in E(M)-x$, $\si(M/z)$ is not $3$-connected with an $N$-minor. Suppose the contrary, by \ref{criterion-1}, $z\in \cl_M(T)$. So, $x$ is in the non-trivial parallel class $\cl_M(T)$ of $M/z$. Since $N$ is simple and $M/z$ has an $N$-minor, then $M/z\del x$, and, therefore, $M\del x$ have an $N$-minor. This is a contradiction to what we proved before. So, $x$ is the unique element of $M$ such that $\si(M/x)$ is $3$-connected with an $N$-minor. By Theorem \ref{teorema-costalonga}, $r(M)-r(N)=1$. 

Consider the structures defined as in Corollary \ref{splitter-plus}. By what we proved, for all choices of $M_1,\dots, M_n$, we have $I=\{x\}$ and $n=3$ or $4$. As $M/x$ has a parallel class with $3$ elements, then $x=x_3$ or $x=x_4$, so, we have two cases to consider:
\smallskip

{\it Case 1. We may pick $M_1,\dots,M_n$ with $x=x_n$: } For all $y\in E(M)-x$ such that $M\del y$ is $3$-connected with an $N$-minor, we have $y\in T$ by \ref{criterion-3}. In particular this holds for each $y\in I\s\cup \cl_M(T)-x$. So, $I\s\cup(\cl_M(T)-x)\cont T$. This implies the validity of the theorem and, in particular, of item (a). Contradicting the assumption that the Theorem fails for $M$, $N$ and $T$.
\smallskip

{\it Case 2. Otherwise: } Now, necessarily, $n=4$ and $x=x_3$. If $M':=M\del x_1,x_2,x_4$ is $3$-connected, then,  as $M'/x$ is $3$-connected, we could choose $x_4=x$, which does not hold in Case 2. Note, that $x$ is in a cocircuit with size at most two in $M'$. As $M_2=M\del x_1,x_2$ is $3$-connected, then $x$ is in a serial pair of $M'$ with an element $z$. This implies that $M\del x_1,x_2,x_4/z\cong N$. By Lemma \ref{wu}, $\si(M/z)$ is $3$-connected with an $N$-minor. This contradicts  the uniqueness of $x$ established before.
\end{claimproof}

\begin{claim}\label{criterion-6}
If $\co(M\del x)$ and $\co(M\del x,y)$ are $3$-connected with $N$-minors, then $\co(M\del y)$ is $3$-connected and $x,y\in T$.
\end{claim}
\begin{claimproof}
Suppose the contrary. By \ref{criterion-3}, $\co(M\del y)$ is not $3$-connected. By the dual of Lemma \ref{w36}, there is a corank-$3$ circuit $C$ containing $x$ with $y\in\cl\s(C)-C$.

First assume that $C\neq T$. If $M$ has a $4$-cocircuit $D\s$ contained in $C\cup y$, then, as $|D\s\cap C|\ge 3$ and $T\ncont C$, there is $z\in (D\s\cap C)-T$. So, $D\s-z$ is a triad of $M\del z$ contained in $\cl\s_{M\del z}(C-z)$ and, by the dual of Lemmas \ref{w37} and \ref{w38}, $\co(M\del z)$ is $3$-connected with an $N$-minor. But this contradicts \ref{criterion-3}. Thus, $C\cup y$ contains no $4$-cocircuit of $M$. But $r\s(C\cup y)=3$ and $y\in \cl\s(C)-C$, so $C\cup y$ is the disjoint union of a singleton set $\{e\}$ and a non-trivial coline $L\s$ containing $y$. By the dual of Lemmas \ref{w37} and \ref{w38} again, $\co(M\del e)$ is $3$-connected with an $N$-minor. By \ref{criterion-3}, $e\in T$. By the dual of Lemma \ref{w37}, for some $f\in L\s$, $\co(M\del x,y)\cong \co(M\del f,y)\cong\co(M\del L\s)$ has an $N$-minor. Thus, $M/y$ has an $N$-minor and, by Bixby's Lemma, $\si(M/y)$ is $3$-connected with an $N$-minor. By \ref{criterion-4}, $y\in T$. Since $T$ meets $L\s$, it follows that $L\s$ is a triad and, as a consequence, $|C|=3$. By orthogonality, there is $g\in (L\s\cap T)-y\cont C$. Since $e\in (C\cap T)-L\s$, then $C\cup T$ is a $4$-segment of $M$ meeting a triad, a contradiction since $M$ is $3$-connected and not isomorphic to $U_{2,4}$. Therefore, $C=T$.

Let $C\s$ be a cocircuit such that $y\in C\s\cont T\cup y$. If $C\s$ is a triad, we have a contradiction to \ref{criterion-2} since $y\in C\s-T$ and $M\del y$ has an $N$-minor. So, $C\s$ is a $4$-cocircuit and $C\s=T\cup y$. By Bixby's Lemma, $\si(M/y)$ is $3$-connected. Since $y\notin T$, then $M/y$ has no $N$-minor by \ref{criterion-4}. If $r(M)=r(N)$, then $N\cong M\del x,y$ and $M\del y$ is $3$-connected, therefore $r(M)-r(N)\in \{1,2\}$. 

For all $2$-subsets $A$ of $T$, $M\del A$ has no $N$-minor because, otherwise, $y$ would be in the serial pair $C\s-A$ of $M\del A$ and $M/y$ would have an $N$-minor.

As $\co(M\del x,y)$ has an $N$-minor, hence $r\s(M)-r\s(N)\ge 2$. If $r\s(M)-r\s(N)\ge 3$, then, by Theorem \ref{wcomut}, there is $z\in E(M)$ such that $\co(M\del z)$ and $\co(M\del x,z)$ are $3$-connected with $N$-minors. By \ref{criterion-3}, $z\in T$. So, for $A:=\{x,z\}\cont T$, $M\del A$ has an $N$-minor, a contradiction. Therefore, $r\s(M)-r\s(N)=2$. To prove the theorem and item (b), we have to find an $N$-minor $N'$ of $M$ with $E(M)-E(N')\cont T\cup y$. Consider a chain of matroids, sets and elements as in Corollary \ref{splitter-plus}. Let $a$ and $b$, in this order, be the elements deleted from $M$ in order to get $M_n$ from $M$ as in the chain (recall that $r\s(M)-r\s(N)=2$). By Lemma \ref{wu}, $\co(M\del a)$ is $3$-connected with an $N$-minor, hence, by \ref{criterion-3}, $a\in T$. It follows from Lemma \ref{wu}, \ref{criterion-5}, \ref{criterion-rankgap} and \ref{criterion-4} that $I\cont T$. We just have to prove now that $b=y$. Suppose the contrary. If $b\in T$, then, for $A:=\{a,b\}\cont T$, $M\del A$ has an $N$-minor, a contradiction, as we saw before. Thus, $b\notin T$ and, by \ref{criterion-3}, $\co(M\del b)$ is not $3$-connected. So, $a$ and $b$ play similar roles as $x$ and $y$ and applying the same steps for $a$ and $b$ as we did for $x$ and $y$, we conclude that $D\s:=T\cup b$ is a cocircuit of $M$. By circuit elimination on $C\s:=T\cup y$, $D\s$ and any element $e$ of $T$, it follows from the cosimplicity of $M$ and from the orthogonality with $T$ that $(T-e)\cup\{y,b\}$ is a cocircuit of $M$. Therefore, $y$ is in a series class of $M\del a,b$, which has an $N$-minor. But this implies that $M/y$ has an $N$-minor, a contradiction. So, $b=y$ and (b) holds.
\end{claimproof}

Now, consider the structures as given by Corollary \ref{splitter-plus}. It follows from \ref{criterion-5}, \ref{criterion-rankgap} and \ref{criterion-4} that $I\cont T$. If $r\s(M)-r\s(N)\le 2$, it follows from \ref{criterion-3} and \ref{criterion-6} that $I\s\cont T$. This implies that $T\cont E(M)-E(N')$ for $N'=M_n$ and the theorem holds in this case. So, we may assume that $r\s(M)-r\s(N)=3$. If $|I|=0$, then $N'=M\del I\s$ and $M\del e$ is $3$-connected with an $N$-minor for all $e\in I\s$ and by \ref{criterion-3}, $I\s\cont T$. So, $|I|\ge 1$. By Lemma \ref{wu} and \ref{criterion-6}, the elements of $I\s$ with the two least indices in $x_1,\dots,x_n$ are in $T$. So, $|T\cap I|\le 1$ and as $I\cont T$, $|I|=1$. Therefore, $n=4$. Since $I\cont T$, $x_1\notin I$ by the simplicity of $M_1$. If $I=\{x_4\}$, then $M\del e$ is $3$-connected for all $e\in I\s$, a contradiction, as before. Therefore, $I=\{x_3\}$ or $\{x_2\}$. This implies that $T=\{x_1,x_2,x_3\}$.

By Theorem \ref{teorema-costalonga}, there is a coindependent $3$-set $J\s$ of $M$ such that $\co(M\del e)$ is $3$-connected with an $N$-minor for all $e\in J\s$. By \ref{criterion-3}, $J\s=T=\{x_1,x_2,x_3\}$. If $T$ meets a triad $T\s$, then, for $f\in T\s-T$ and $e\in T\cap T\s$ we have that $M\del e$ and, therefore, $M\del e/f$, have an $N$-minor. But in this case, by Lemma \ref{w38}, $\si(M/f)$ is $3$-connected with an $N$-minor, a contradiction to \ref{criterion-5}. Thus, $T$ meets no triads of $M$. Next, we check:

\begin{claim}\label{criterion-7}
We may not pick $M_1,M_2,M_3,M_4$ in such a way that $I=\{x_3\}$.
\end{claim}
\begin{claimproof}
Suppose the contrary. Then, we may pick the chain of matroids in such a way that $M_4=M\del x_1,x_2/x_3\del x_4$ with $I=\{x_3\}$. As $x_4\notin T$, by \ref{criterion-3},  $\co(M\del x_4)$ is not $3$-connected and $M\s$ has a vertical $3$-separation $(A,x_4,B)$. This is, both $A$ and $B$ are $3$-separating sets of $M$, $x_4\in \cl\s(A)\cap\cl\s(B)$ and $r\s(A),r\s(B)\ge 3$. So, $(A,B)$ is a $2$-separation of $M\del x_4$, but $M\del x_4,x_1,x_2/x_3$ is $3$-connected and we may assume, therefore, that $|A-\{x_1,x_2,x_3\}|=|A-T|\le 1$. As $|A|\ge 3$, then $|A\cap T|\ge 2$ and $A$ spans $T$. This implies that $Y=A\u T$ is a $3$-separating set of $M$. Moreover, $Y=T$ or $Y=T\cup y$ for some $y\in E(M)-T$. 

Next we prove that $T\cup x_4$ is a cocircuit of $M$. Note that $x_4\in\cl\s(Y)$. If $Y=T$, then, as $T$ meets no triads, it follows that $T\cup x_4$ is a $4$-cocircuit. So, we may assume that $Y=T\cup y$ for some $y\in E(M)-T$. If $r_M(Y)=2$, $Y$ is a $4$-segment of $M$. As $M/x_3$ has an $N$-minor, then so has $M\del y$. In this case, $M\del y$ is $3$-connected, contradicting \ref{criterion-3} since $y\notin T$. Thus, $r_M(Y)=3$. Since $|Y|=4$ and $Y$ is $3$-separating, it follows that $r\s(Y)=3=r\s(T)$. Now, $T$ cospans $y$, and as $T\cup y$ cospans $x_4$, it follows that $T$ cospans $x_4$ and $T\cup x_4$ is a cocircuit of $M$ since $T$ meets no triads.

As $T\cup x_4$ is a cocircuit of $M$, hence $T-\{x_3,x_4\}$ is a serial pair of $M_2=M\del x_1,x_2$ which is $3$-connected with at least $4$ elements, a contradiction. Thus, \ref{criterion-7} holds.
\end{claimproof}

Now, By \ref{criterion-7}, $I=\{x_2\}$ for all choices of chains. This implies that there are no pair of elements $\{a,b\}\cont E(M)$ such that $M\del a$ and $M\del a,b$ are $3$-connected with an $N$-minor. In particular, $M\del x_1,x_3$ is not $3$-connected. But $M_3=M\del x_1,x_3/x_2$ is $3$-connected. As $M\del x_1$ is $3$-connected, then $x_2$ is in a serial pair $\{x_2,z\}$ of $M\del x_1,x_3$. Hence, $M/z$ has an $N$-minor. Since $M\del x_1$ is $3$-connected, it follows that $\{z,x_2,x_3\}$ is a triad of $M\del x_1$. But $T$ meets no triads of $M$ and, therefore, $C\s:=\{z,x_1,x_2,x_3\}=T\cup z$ is a $4$-cocircuit of $M$. If $z\in\cl(T)$, then $r\s(C\s)=2$ and $C\s$ is a $2$-separating set of $M$. This implies that $r(M)=r(C\s)=2$, contradicting the fact that $r(M)=r(N)+1\ge 3$. So, $z\notin \cl(T)$ and $T$ is a triangle of $M/z$ contained in $C\s$. By Lemma \ref{w38}, $\si(M/z)$ is $3$-connected. But $M/z$ has an $N$-minor. By \ref{criterion-4}, $z\in T$, a contradiction. This proves the theorem.
\end{proofof}

\begin{proofof}{\emph Proof of Theorem \ref{binary-criterion}: } If $|E(N)|\ge 4$, then, since binary matroids cannot contains $4$-point lines nor circuits meeting cocircuits in an odd number of elements, then items (a) and (b) of Theorem \ref{criterion} does not hold. In this case, Theorem \ref{criterion} implies the theorem. So we may assume that $|E(N)|\le 3$. If $N\ncong U_{1,3}$, then, since $N$ is $3$-connected, if follows that $N$ is isomorphic to a minor of $M':=M|T\cong U_{2,3}$. So, we may assume that $N\cong U_{1,3}$. In this case $M$ has a circuit $C$ meeting $T$ but different from $T$. Now, for any $e\in C-T$, the theorem holds for $M'=(M|T\cup C)/(C-(T\cup e))$.\end{proofof}

We define $K_{3,3}^{1,1}$ as the graph in Figure \ref{fig-K}. The following lemma is a well-known result and is a straightforward consequence of Seymour's Splitter Theorem.

\begin{lemma}\label{k5minor}
If $G$ is a $3$-connected graph with a $K_5$-minor then, either $G\cong K_5$ or $G$ has a $K_{3,3}^{1,1}$-minor.
\end{lemma}

\begin{center}
 \begin{figure}[h]\centering
  \begin{tikzpicture}
   \tikzstyle{node_style} =[shape = circle,fill = black,minimum size = 2pt,inner sep=1pt]
   \node[node_style] (d) at (0,0) {};     \node () at (-0.2,0) {$d$};
   \node[node_style] (c) at (0,1) {};     \node () at (-0.2,1) {$c$};
   \node[node_style] (u) at (0,2) {};     \node () at (-0.2,2) {$u$};
   \node[node_style] (b) at (1,0) {};     \node () at (1.2,0) {$b$};
   \node[node_style] (a) at (1,1) {};     \node () at (1.2,1) {$a$};
   \node[node_style] (v) at (1,2) {};     \node () at (1.2,2) {$v$};
   \draw (a)--(u)--(b)--(a)--(c)--(v)--(d)--(c)--(b)--(d)--(a)--(u)--(v);
  \end{tikzpicture}
  \caption{$K\cong K_{3,3}^{1,1}$}\label{fig-K}
 \end{figure}
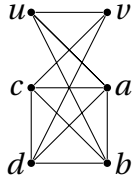
\end{center}


\begin{proofof}{\it Proof of Theorem \ref{k_5 rounded}: }
We have to prove that for each $3$-connected simple graph $G$ with a $K_5$-minor and for each triangle $T$ of $G$, $G$ has a $K_5$-minor using $E(T)$. Consider a counter-example $G$ with $|E(G)|$ as small as possible. By Theorem \ref{binary-criterion}, we may assume that $E(G)-E(K_5)\cont E(T)$. As no edges may be added to $K_5$ in order to get a $3$-connected simple graph, then $|V(G)|=6$ or $7$. 

First suppose that $|V(G)|=7$. Now, $G$ is obtained from $K_5$ by two vertex splittings by edges $e_1$ and $e_2$ and, possibly, by the addition of an edge $e_3$ with $e_1,e_2\in E(T)$ and $e_3\in E(T)$ in case it is added. So, $e_1$ and $e_2$ are adjacent in $G$. When obtaining $G$ in this process, we split a vertex by $e_1$, then we split an endvertex of $e_1$ by $e_2$. If $x$, $y$, and $z$ are the three vertices incident to $e_1$ or $e_2$ in after these two splittings, then the sum of their degrees is $8$. So, in order to get $G$ with all vertices with degree at least $3$, the addition of $e_3$ is necessary. Note that this process describes that $G$ is obtained from $K_5$ by expanding a vertex into the triangle $T$ with vertices $x$, $y$ and $z$. Observe that there are two vertices $u,v\in V(T)$ with degree $3$. It is clear that for the edges $e,f\notin E(T)$ incident to $u$ and $v$ respectively, we have $\si(G/e,f)\cong K_5$. 

So, we may assume that $|V(G)|=6$. By Lemma \ref{k5minor}, up to labels, $G$ is obtained from $K\cong K_{3,3}^{1,1}$ (the graph in Figure \ref{fig-K}), by adding the edges of $E(T)-E(K)$. Since $K/uv\cong K_5$, then $uv\in T$ and we may assume without losing generality that $V(T)=\{u,v,a\}$, so $G=K+va$. Now, it is clear that $G\del ba/ub$ is a $K_5$-minor of $G$ using $T$. This proves the theorem.\end{proofof}

The following lemma has a slightly stronger conclusion than \cite[Proposition 9.3.5]{Oxley} (beyond the conclusions of \cite[Proposition 9.3.5]{Oxley}, it states beyond that $R$ has a $K$-minor and describes the way it is obtained), but the proof for \cite[Proposition 9.3.5]{Oxley} also holds for the following lemma.

\begin{lemma}\label{oxley935}
Let $R=K\oplus_3 L$ be a $3$-sum of binary matroids, where $K$ and $L$ are $3$-connected and $E(K)\cap E(L)=S$. Then there are $X,Y\cont E(L)-S$ such that $K\cong R/X\del Y$, where $K$ is obtained from $R/X\del Y$ by relabeling the elements $s_1$, $s_2$ and $s_3$ of $S$ in $K$ by respective elements $l_1$, $l_2$ and $l_3$ of $L$.
\end{lemma}

\begin{proofof}{\it Proof of Theorem \ref{regular}: } Suppose that $R$ is a matroid contradicting the theorem minimizing $|E(R)|$. So, $R$ has a triangle $T$ and an $M$-minor for some $M\in \F$ but $M$ has no $\F$-minor using $T$. If $R$ is graphic or cographic, the theorem holds for $R$, so assume the contrary. By Seymour's Decomposition Theorem for Regular Matroids~\cite[Theorem 13.1.1]{Oxley}, there are matroids $K$ and $L$ with at least $7$ elements each, intersecting in a common triangle $S$ such that $R=K\oplus_3 L$ with $L$ being $3$-connected and $K$ being $3$-connected up to parallel classes of size two meeting $S$. Under these circunstances, we may assume that $|E(K)\cap E(M)| \ge |E(L)\cap E(M)|$.

If $C$ is a cycle of $R$ meeting both $E(K)$ and $E(L)$, then there is $s\in S$ such that $(C\cap E(N))\cup s$ is a cycle of $N$ for $N=K,L$. As we picked $L$ with no parallel pairs, it follows that $\cl_R(E(K)-S)\cap E(L)=\emptyset$.

Let us first check that $K$ has an $M$-minor. Let $M=R/I\del I\s$ for some independent set $I$ and coindependent set $I\s$ of $R$. Since $\lambda_M(E(K)\cap E(M))\le \lambda_R(E(K)\cap E(M))=2$, then as $M$ is internally $4$-connected, it follows that $|E(M)\cap E(L)|\le 3$, and, moreover, $E(M)\cap E(L)$ is not a triad of $M$ because $M$ has no triads. This implies that $E(M)\cap E(L)\cont \cl_M(E(M)-E(L))$. By the format of the family of circuits of $R$, it follows that $E(M)\cap E(L)\cont \cl_R(E(K)-S)$, which is empty. So, $E(M)\cont E(K)$.
By Lemma \ref{oxley935}, there is a minor $K'$ of $R$ obtained by relabeling the elements $s_1$, $s_2$ and $s_3$ of $S$ in $K$ by respective elements $l_1$, $l_2$ and $l_3$ of $L$. Consider the matroid $K''$ obtained from $K'$ by contracting each $l_i$ for those indices $i\in\{1,2,3\}$ such that $s_i\in \cl_L(I\cap E(L))$. Now $K''$ is obtained from $R/(E(L)\cap I)\del (E(L)\cap I\s)$ by relabeling the remaining elements of $S$. This implies that $K''$ and, therefore, $K$ have $M$-minors.

If $T\cont E(K)$, then $K$ has an $\F$-minor using $T$ by the minimality of $R$. But $R$ has an $K$-minor using $T$ by Lemma \ref{oxley935} and this implies the theorem. So, $T$ meets $E(L)$. As $\cl_R(E(K)-S)\cap E(L)=\emptyset$, it follows that $X:=T\cap E(L)$ has at least two elements. As $L$ is $3$-connected and $\lambda_L(S)=\lambda_L(X)=2$, then $\kappa_L(S,X)=2$. By Tutte's Linking Theorem~\cite[Theorem 8.5.2]{Oxley}, there is a minor $N$ of $L$ with $E(N)=S\cup X$ such that $\lambda_N(S)=2$. Hence: 
\[2=\lambda_N(S)=r_N(S)+r_N(X)-r(N)\le 4-r(N).\]

So, $r(N)\le 2$. But $r(N)\ge r_N(S)\ge \lambda_N(S)=2$. Also $r_N(X)\ge \lambda_N(X)=2$. This implies that $S$ spans $N$ and $X$ contains no parallel pairs of $N$. Now, each element of $X$ is in parallel with an element of $S$ in $N$. Therefore, for $N=L/A\del B$, we have that $R/A\del B$ is obtained from $K$ by relabeling the elements of $S$ by elements of $T$. So, $R/A\del B$ is $3$-connected with $T$ as triangle and has an $M$-minor. By the minimality of $R$, $R/A\del B$ has an $\F$-minor using $T$ and this proves the lemma.
\end{proofof}

\section{Sharpness}\label{sec-sharp}

In this section we establish that Theorem \ref{criterion} is sharp in the sense that all described cases may occur indeed. 

First we construct an example for Theorem \ref{criterion} with $E(M)\cont E(N)\cup T$. Consider a complete graph $K$ on $n\ge 14$ vertices. Let $X:=\{v_{i,j}:i=1,2,3$ and $ j=1,2,3,4\}$ be a $12$-subset of $V(K)$. Consider a triangle $T$ on vertices $u_1$, $u_2$ and $u_3$, disjoint from $K$. Let $G=K\cup T+\{u_iv_{i,j}:i=1,2,3$ and $j=1,2,3,4\}$. Define, for disjoint subsets $A$ and $B$ of $E(T)$, $H:=G\del A/B$, $M:=M(G)$ and $N:=M(H)$ provided $H$ is a simple graph. For each $e\in E(G)-E(T)$, $G/e$ has at least $3$ parallel pairs; so, $|E(\si(M/e))|\le |E(G)|-4<|E(H)|$ and, therefore $M/e$ has no $N$-minor. For $e$ incident to $v\in V(K)$, in order to get a minor of $G\del e$ with $|V(K)|-12$ vertices with degree $|V(K)|-1\ge 13$ and $12$ vertices with degree $|V(K)|\ge 14$ as in $H$, it is necessary to contract some edge out of $T$, thus $G\del e$ has no $H$-minor either. Note that we may pick $A$ with any size from $0$ to $3$, We always may pick $B=\emptyset$, and provided $|A|\ge 1$ we may pick $B$ with size from $1$ to $3-|A|$.

Now, let us construct an example satisfying item (a) of Theorem \ref{criterion}. Let us pick $M$ as a restriction of the affine space $\mathbb{R}^3$. Consider a $4$-subset $L:=\{a,b,c,x\}$ of an line $R$. Let $T:=L-x$. Now consider for each $y\in T$ a line $R_y$ meeting $L$ in $y$ in such a way that no three lines among $R$, $R_a$, $R_b$ and $R_c$ lay in a same plane. Let $m\ge 6$. For each $y\in T$, pick a $m$-subset $L_y$ of $R_y$ containing $y$. Let $M$ be the restriction of the affine space to $L\cup L_a\cup L_b\cup L_c$. Let $N=M/x\del a,b$ or $N=M/x\del T$. Note that it is not possible to get a rank-$3$ minor of $M$ with $3$ disjoint $(m-1)$-segments by contracting an element other than $x$. So all $N$-minors of $M$ are minors of $M/x$ and therefore, deleting at least two elements of $T$ is necessary to get an $N$-minor since $T$ is a parallel class of $M/x$. So, this is the unique way to get an $M/x\del a,b$-minor of $M$. Moreover, deleting the element in the intersection of the three $m$-lines is the unique way to get and $M/x\del T$ from an $M/x\del a,b$-minor of $M$.

Next, we construct an example satisfying item (b) of Theorem \ref{criterion}. We denote by $M+e$ the matroid obtained by adding $e$ freely to $M$. Start with a projective geometry $P$ with $r(P)\ge 6$. Let $F$ be a flat of $P$ with $4\le r(F)\le r(P)-2$. Consider a copy $U$ of $U_{2,4}$ on ground set $T\cup x:=\{x,x_1,x_2,x_3\}$ with $(T\cup x)\cap E(P)=\emptyset$. Let $y$ be an element out of $E(P)\cup T\cup x$. Let $M$ be the matroid obtained by adding $y$ freely to the flat $F\cup T$ of $(P+x)\oplus_2 U$. Note that $E(P)$ is a hyperplane of $M$ and, therefore $T\cup y$ is a $4$-cocircuit of $M$. Define $N_1=M\del y/x_2\del x_3$ and $N_2:=N_1/x_1$. Note that $N_1=P+x_1$ and $N_2$ is the truncation of $P$ with rank $r(P)-1$ (the matroid on $E(P)$ whose independent sets are those independent sets of $P$ with size at most $r(P)-1$).

Let $i\in\{1,2\}$ and $N=M/X\del Y$ be an $N_i$-minor of $M$. Note that $r\s(M)-r\s(N)=2$ and $r(M)-r(N)\in\{1,2\}$. For each $p\in E(P)$, $|E(si(M/p))|<|E(N_i)|$. Thus, no element $p\in P$ may be contracted in $M$ in order to get an $N_i$-minor. So, $X\cont T\cup y$.

Let us check that $T\cup y$ meets no circuit of $M$ with less than six elements other than $T$. Indeed, $M\del y$ is a two sum of a $4$-point line on $T\cup x$ and a matroid with rank greater than five with $x$ as free element. Thus all circuits of $M\del y$ meeting $T$, except for $T$ itself, have more than five elements. Moreover, $M$ is obtained from $M\del y$ adding $y$ as a free element to a flat with rank greater than $4$ and, therefore, all circuits of $M$ containing $y$ also have more than five elements. Hence, the triangles of $M/X$ are precisely the triangles of $P$. Moreover, those must be the same triangles of $N$ since all triangles of $N$ are triangles of $M/X$ and they occur in the same number. As deleting an element of $P$ from $M/X$ would result in a matroid with less triangles than $N_i$, it follows that $Y\cont T\cup y$. Hence, $E(P)\cont E(N)$ for each minor $N$ of $M$ isomorphic to $N_1$ or $N_2$. 

Let us check that $M/y$ has no $N_2$-minor and, therefore, no $N_1$-minor too. Suppose for a contradiction that $N$ is an $N_2$-minor of $M/y$. We may assume that $N=M/y,x_1\del x_2,x_3$. Note that $x_1$ is a free element of the rank-$r_P(F)$ flat $F\cup x_1$ of $M/y\del x_2, x_3$. This implies that $N|F$ is a truncation of rank $r_P(F)-1$ of the rank-$r_P(F)$ projective geometry $F$. But, as $N_2$ is the rank-$(r(P)-1)$ truncation of $P$ and $r_P(F)\le r(P)-2$, then all rank-$(r_P(F)-1)$ flats of $N_1$ are projective geometries and so is $F$, a contradiction. 

Now, for $i=1,2$, each $N_i$-minor of $M$ is in the form $M\del y,x_i/A$ with $1\le i\le 3$ and $A$ being an $i$-subset of $T-x_i$. Let $A$ be a $2$-subset of $T$. As we proved for $A=\{x_2,x_3\}$, it follows that $M\del A$ has no $N_2$-minor. Moreover, for $x_k\in T-A$ it is clear that $y$ is not a free element of $M\del A/x_k$, which, therefore, is not isomorphic to $N_1$. Thus, $M\del A$ has no minor isomorphic to $N_1$. This implies that $M\del A$ has no minor isomorphic to $N_2$ neither.

\section{Acknowledgements}
The authors thank the Visiting Scholar Fund of the Department of Mathematics and Statistics at Wright State University for the support.


\end{document}